\documentclass[12pt]{amsart}
\usepackage{latexsym}
\usepackage{amsthm}
\usepackage{amsmath}
\usepackage{amsfonts}
\usepackage{amssymb}
\usepackage[dvips]{graphicx}
\usepackage{xypic}
\input xy
\xyoption{all}

\newtheorem{theo}{Theorem}[section]
\newtheorem{lemma}[theo]{Lemma}
\newtheorem{propo}[theo]{Proposition}
\newtheorem{defi}[theo]{Definition}
\newtheorem{coro}[theo]{Corollary}
\newtheorem{rem}[theo]{Remark}
\newtheorem{exam}[theo]{Example}

\newcommand\Ab{\operatorname{\bf Ab}}
\newcommand\Gr{\operatorname{\bf Gr}}
\newcommand\Lin{\operatorname{\bf Lin}}

\newcommand\Ind{\operatorname{Ind}}

\newcommand\Str{\operatorname{\bf Str}}

\newcommand\rank{\operatorname{rank}}

\newcommand\op{\operatorname{op}}

\newcommand\id{\operatorname{id}}
\newcommand\Id{\operatorname{Id}}

\newcommand\Mod{\operatorname{\bf Mod}}
\newcommand\Elem{\operatorname{\bf Elem}}
\newcommand\Emb{\operatorname{\bf Emb}}
\newcommand\Set{\operatorname{\bf Set}}

\newcommand\Pos{\operatorname{\bf Pos}}

\newcommand\ca{\mathcal {A}}

\newcommand\cc{\mathcal {C}}
\newcommand\cd{\mathcal {D}}

\newcommand\ck{\mathcal {K}}
\newcommand\cl{\mathcal {L}}

\newcommand\cw{\mathcal {W}}

\date{November 12, 2011}
 
\begin{document}
\title{Abstract elementary classes and accessible categories}
\author[T. Beke and J. Rosick\'{y}]
{T. Beke$^*$ and J. Rosick\'{y}$^{**}$}
\thanks{$^*$ Supported by the NSA under Grant H98230-11-1-0172, $^{**}$ Supported by MSM 0021622409.} 
\address{\newline T. Beke\newline
Department of Mathematics\newline
University of Massachusetts Lowell\newline
One University Avenue, Lowell, MA 01854, U.S.A.\newline
tibor\_beke@uml.edu
\newline J. Rosick\'{y}\newline
Department of Mathematics and Statistics\newline
Masaryk University, Faculty of Sciences\newline
Kotl\'{a}\v{r}sk\'{a} 2, 611 37 Brno, Czech Republic\newline
rosicky@math.muni.cz
}
\begin{abstract}
We investigate properties of accessible categories with directed colimits and their relationship with categories arising from Shelah's Abstract 
Elementary Classes.  We also investigate ranks of objects in accessible categories, and the effect of accessible functors on ranks.  
\end{abstract}

\maketitle

\section{Introduction}
M. Makkai and R. Par\' e \cite{MP} introduced accessible categories as categories sharing two typical properties of categories
of structures described by infinitary first-order theories -- the existence of sufficiently many directed colimits and the existence 
of a set of objects generating all objects by means of distinguished colimits. Their, purely category-theoretical, definition 
has since then found applications in various branches of mathematics. Very often, accessible categories have all directed colimits. These arise 
as categories of models and elementary embeddings of infinitary theories with finitary quantifiers. In model theory, S.\ Shelah went in a similar 
direction and introduced abstract elementary classes as a formalization of properties of models of generalized logics with finitary quantifiers. 
Our aim is to relate these two approaches. In Section 5, we introduce a hierarchy of accessible categories with directed colimits. The main result 
of that section, Corollary \ref{cor5.7}, sandwiches Shelah's Abstract Elementary Classes between two natural families of accessible categories.

Unlike abstract elementary classes, accessible categories are not equi\-pped with canonical `underlying sets'. Nonetheless, there exists
a good substitute for `size of the model', namely, the \emph{presentability rank} of an object, that can be expressed purely 
in the language of category theory, i.e.\ in terms of objects and morphisms.  In this sense, accessible categories take an extreme 
`signature-free' and `elements-free' view of abstract elementary classes.  From this point of view, Shelah's Categoricity 
Conjecture, the driving force of abstract elementary classes (see \cite{B}), turns on the subtle interaction between ranks of objects 
and directed colimits in accessible categories.

The connection between abstract elementary classes and accessible categories was discovered, independently, by M.\ Lieberman as well; 
see \cite{L} and \cite{L1}. Our Corollary \ref{cor5.7} simplifies Lieberman's description, \cite{L1}, Proposition 4.8 and Claim 4.9.

\section {Accessible categories}
In order to define accessible categories one just needs the concept of a $\lambda$-directed colimit where
$\lambda$ is a regular cardinal number. This is a colimit over a diagram $D:\cd\to\ck$ where $\cd$ is a
$\lambda$-directed poset, considered as a category. An object $K$ of a category $\ck$ is called 
$\lambda$-\textit{presentable} if its hom-functor $\hom(K,-):\ck\to\Set$ preserves $\lambda$-directed
colimits; here $\Set$ is the category of sets.

A category $\ck$ is called $\lambda$-\textit{accessible}, where $\lambda$ is a regular cardinal, provided that
\begin{enumerate}
\item[(1)] $\ck$ has $\lambda$-directed colimits,
\item[(2)] $\ck$ has a set $\ca$ of $\lambda$-presentable objects such that every object of $\ck$ is a 
$\lambda$-directed colimit of objects from $\ca$.
\end{enumerate}
A category is \textit{accessible} if it is $\lambda$-accessible for some regular cardinal $\lambda$.  

A \textit{signature} $\Sigma$ is a set of (infinitary) operation and relation symbols. These symbols are
$S$-sorted where $S$ is a set of sorts. It is advantageous to work with many-sorted signatures but it is
easy to reduce them to single-sorted ones. One just replaces sorts by unary relation symbols and adds axioms
saying that there are disjoint and cover the underlying set of a model. Thus the underlying set of an $S$-sorted
structure $A$ is the disjoint union of underlying sets $A_s$ over all sorts $s\in S$. $|A|$ will denote the cardinality 
of the underlying set of the $\Sigma$-structure $A$. The category of all $\Sigma$-structures and homomorphisms (i.e., mappings preserving 
all operations and relations) is denoted by $\Str(\Sigma)$. A homomorphism is called a substructure embedding if it is injective 
and reflects all relations. Any inclusion of a substructure is a substructure embedding. Conversely, if $h:A\to B$ is a substructure embedding then $A$ 
is isomorphic to the substructure $h(A)$ of $B$. The category of all $\Sigma$-structures and substructure embeddings is denoted by $\Emb(\Sigma)$. 
Both $\Str(\Sigma)$ and $\Emb(\Sigma)$ are accessible categories, cf.\ \cite{AR}, 5.30 and 1.70.  

A signature $\Sigma$ is \textit{finitary} if all relation and function symbols are finitary. For a finitary signature,
the category $\Str(\Sigma)$ is locally finitely presentable and $\Emb(\Sigma)$ is finitely accessible. In both cases, 
there is a cardinal $\kappa$  such that, for each regular cardinal $\kappa\leq\mu$, a $\Sigma$-structure $K$ is 
$\mu$-presentable if and only if $|K|<\mu$. This follows from the downward L\" owenheim-Skolem theorem; see \cite{MP}, 3.3.1.

Let $\lambda$ be a cardinal and $\Sigma$ be a $\lambda$-ary signature, i.e., all relation and function symbols have arity 
smaller than $\lambda$. Given a cardinal $\kappa$, the language $L_{\kappa\lambda}(\Sigma)$ allows less than $\kappa$-ary 
conjunctions and disjunctions and less than $\lambda$-ary quantifications.  A substructure embedding of $\Sigma$-structures is called
$L_{\kappa\lambda}$-\textit{elementary} if it preserves all $L_{\kappa\lambda}$-formulas. A \textit{theory} $T$ is a set 
of sentences of $L_{\kappa\lambda}(\Sigma)$. $\Mod(T)$ denotes the category of $T$-models and homomorphisms, $\Emb(T)$ 
the category of $T$-models and substructure embeddings while $\Elem(T)$ will denote the category of $T$-models and $L_{\kappa\lambda}$-elementary 
embeddings. The category $\Elem(T)$ is accessible (see \cite{AR}, 5.42).  For certain theories $T$, the category $\Mod(T)$ 
does not have $\mu$-directed colimits for any regular cardinal $\mu$ and thus fails to be accessible.

A theory $T$ is called \textit{basic} if it consists of sentences
$$
(\forall x)(\varphi(x)\Rightarrow\psi(x))
$$
where $\varphi$ and $\psi$ are positive-existential formulas and $x$ is a string of variables.  For a basic theory $T$, 
the category $\Mod(T)$ is accessible. Conversely, every accessible category is equivalent to the category of models 
and homomorphisms of a basic theory. All these facts can be found in \cite{MP} or \cite{AR}.  
 
\textit{Locally presentable categories} are defined as cocomplete accessible categories. Following \cite{AR}, 1.20, 
each locally $\lambda$-presentable category is locally $\mu$-presentable for each regular cardinal $\mu>\lambda$. 
Let $\lambda$ be an uncountable regular cardinal. The category $\Pos_\lambda$ of $\lambda$-directed posets 
and substructure embeddings is $\lambda$-accessible but it is not $\mu$-accessible for all regular cardinals $\mu>\lambda$.  
Following \cite{MP} 2.3, let us write $\lambda\triangleleft\mu$ whenever $\Pos_\lambda$ 
is $\mu$-accessible.   There are arbitrarily large regular cardinals $\mu$ such that $\lambda\triangleleft\mu$ 
and, at the same time, arbitrarily large regular cardinals $\mu$  such that $\lambda\triangleleft\mu$ does not hold. Thus the \textit{accessibility 
spectrum} of $\Pos_\lambda$ has a proper class of gaps. Generally, if $\lambda\triangleleft\mu$ then any $\lambda$-accessible category $\ck$ 
is $\mu$-accessible; see \cite{MP}, 2.3.10 or \cite{AR}, 2.11.  By \cite{AR}, 2.13 (1), one has $\omega\triangleleft\lambda$ for every uncountable 
regular cardinal $\lambda$. Thus a finitely accessible category is $\mu$-accessible for all uncountable regular cardinals $\mu$.
 
\begin{defi}\label{def2.1}
{
\em
We say that a category $\ck$ is \textit{well} $\lambda$-\textit{accessible} if it is $\mu$-accessible for each regular
cardinal $\lambda\leq\mu$. 

$\ck$ is \textit{well accessible} if it is well $\lambda$-accessible for some regular cardinal $\lambda$.
}
\end{defi}

We have just seen that any locally presentable category and any finitely accessible category is well-accessible.

\begin{defi}\label{def2.2}
{
\em
Let $\lambda$ be a regular cardinal.  We say that an object $K$ of a category $\ck$ has \textit{presentability rank} (or, for brevity, \textit{rank}) 
$\lambda$ if it is $\lambda$-presentable but not $\mu$-presentable for any regular cardinal $\mu<\lambda$.  We will write $\rank(K)=\lambda$.
}
\end{defi}

This concept was introduced by Makkai and Par\'e  under the name of \emph{presentability} of an object.  See the very last line of p.~29 of \cite{MP}.  
 
\begin{rem}\label{re2.3}
{
\em
(1) We have mentioned that in $\Emb(\Sigma)$, $\Sigma$ finitary, there is a cardinal $\kappa$ such that, for each cardinal 
$\kappa\leq\mu$, a $\Sigma$-structure $A$ is $\mu$-presentable if and only if $|A|<\mu$.  Therefore for all large enough 
structures, $\rank(A)=|A|^+$, i.e., presentability ranks play the role of cardinalities.

(2) Let $T$ be a basic theory in a language $L_{\kappa\omega}(\Sigma)$. Since $\Emb(T)$ is closed under directed colimits in $\Emb(\Sigma)$,
$\rank(X)\leq|X|^+$ for each $T$-model $X$.  By the downward L\" owenheim-Skolem theorem, cf.\ \cite{MP} 3.3.1, there is a cardinal $\mu$
such that, for $\mu\leq|X|$, we have $\rank(X)=|X|^+$.  Moreover, either the presentability ranks of $T$-models in $\Emb(T)$ form a set or $T$ has models
of presentability rank $\nu^+$ for all $\mu\leq\nu$. This means that the \textit{presentability spectrum} of $\Emb(T)$ does not have arbitrarily large 
gaps.

(3) If $T$ is a basic theory in a language $L_{\kappa\omega}(\Sigma)$ then $\Mod(T)$ does not need to be closed under directed colimits
in $\Str(\Sigma)$. But this holds for basic theories $T$ in a language $L_{\kappa\omega}^{\ast}(\Sigma)$ allowing only finitary conjuctions.
Thus $\Mod(T)$ is an accessible category with directed colimits in this case. \cite{MP} calls categories equivalent to such $\Mod(T)$
$\infty,\omega$-\textit{elementary}. Every finitely accessible category is $\infty,\omega$-\textit{elementary}; see \cite{MP} 4.3.2.
The downward L\" owenheim-Skolem theorem applies to $\infty,\omega$-elementary $\Mod(T)$, and implies, similarly to (2), that $\rank(X)$ eventually 
coincides with $|X|^+$ and the presentability spectrum of $\Mod(T)$ does not have arbitrarily large gaps.

Even if $T$ is only a basic theory in a language $L_{\kappa\omega}(\Sigma)$, then $\Mod(T)$ is closed under $\lambda$-directed colimits
in $\Str(\Sigma)$ for some regular cardinal $\lambda$. Then the downward L\" owenheim-Skolem theorem still implies that $\rank(X)$ eventually 
coincides with $|X|^+$ and the presentability spectrum of $\Mod(T)$ does not have arbitrarily large gaps.

(4) Let $T$ be an arbitrary theory in a language $L_{\kappa\omega}(\Sigma)$. We know that $\Elem(T)$ is an accessible category and,
following \cite{AR} 5.39, $\Elem(T)$ is closed under directed colimits in $\Str(\Sigma)$. Thus $\Elem(T)$ is an accessible category with 
directed colimits. Moreover, $\rank(X)=|X|^+$ in $\Elem(T)$ for all sufficiently large $T$-models $X$ and the presentability spectrum 
of $\Elem(T)$ does not have arbitrarily large gaps.
 
Following \cite{MP} 3.2.8, there is a basic theory $T'$ in another language $L_{\kappa\omega}^{\ast}(\Sigma')$ such that 
the categories $\Elem(T)$ and $\Mod(T')$ are equivalent.  Thus $\Elem(T)$ is $\infty,\omega$-elementary.

(5) Assuming GCH, in any accessible category $\ck$, $\rank(K)$ is a successor cardinal for all objects $K$ (with the possible 
exception of a set of isomorphism types).  Indeed, let $\ck$ be a $\kappa$-accessible category and $\rank(K)=\lambda$ 
with $\kappa<\lambda$. $\lambda$ is a regular cardinal by definition; if $\lambda$ was a limit cardinal then (as it is 
uncountable)\ it would be a weakly inaccessible cardinal. Since GCH is assumed, $\lambda$ is inaccessible. Thus, given 
$\alpha < \kappa <\lambda$ and $\beta <\lambda$, we have $\beta^\alpha <\lambda$. Following \cite{AR} 2.13 (4), 
$\kappa\triangleleft\lambda$.  Following \cite{MP}, 2.3.11, the object $K$ can be exhibited as a $\kappa$-directed colimit 
of $\kappa$-presentable objects along a diagram of size less than $\lambda$.  Let the size of that diagram be $\nu$; then 
$K$ is $\nu^+$-presentable. That would mean $\nu^+=\lambda$, contradicting that $\lambda$ is a limit cardinal.
}
\end{rem}

By `category' we always mean a locally small one, i.e., having a set of morphisms between any two objects. A category is 
called \textit{small} if it has a set of objects. We say that a category is \textit{large} if it is not equivalent to a small category.
This means that it has a proper class of non-isomorphic objects.  In Remark \ref{re2.3} we saw several families of large accessible categories that, 
starting from some cardinal, possess objects of every possible presentation rank. Since in those examples, this property followed 
from the downward L\" owenheim-Skolem theorem, we will call such categories LS-accessible. 

\begin{defi}\label{def2.4}
{
\em
An accessible category $\ck$ will be called $\lambda$-\textit{LS-accessible} if $\ck$ has an object of presentability rank 
$\mu^+$ for each cardinal $\mu\geq\lambda$.

$\ck$ is \textit{LS-accessible} if it is $\lambda$-LS-accessible for some cardinal $\lambda$.
}
\end{defi} 

We will deal with this concept later. However, we have not been able to find any large accessible category which is not LS-accessible.
 
\section {Accessible functors}

A functor $F:\ck\to\cl$ is called $\lambda$-accessible if $\ck$ and $\cl$ are $\lambda$-accessible categories and $F$ preserves 
$\lambda$-directed colimits. $F$ is called accessible if it is $\lambda$-accessible for some regular cardinal $\lambda$. 
By the Uniformization theorem (see \cite{MP}, 2.5.1 or \cite{AR}, 2.19), for each accessible functor $F$ there is a regular 
cardinal $\lambda$ such that $F$ is $\lambda$-accessible and \textit{preserves} $\lambda$-\textit{presentable objects}; it means 
that if $K$ is $\lambda$-presentable then $F(K)$ is $\lambda$-presentable. The same is then true for each $\mu\triangleright\lambda$.  

\begin{defi}\label{def3.1}
{
\em
A functor $F:\ck\to\cl$ will be called \textit{well} $\lambda$-\textit{accessible} if $\ck$ and $\cl$ are well $\lambda$-accessible 
categories and $F$ preserves $\mu$-directed colimits and $\mu$-presentable objects for each $\lambda\leq\mu$. $F$ is called 
\textit{well acce\-ssib\-le} if it is well $\lambda$-accessible for some regular cardinal $\lambda$. 
}
\end{defi} 

\begin{rem}\label{re3.2}
{
\em
(1) Of course, whenever $F$ preserves $\lambda$-directed colimits then it preserves $\mu$-directed ones for $\mu\geq\lambda$.

(2) Every colimit preserving functor between locally $\lambda$-presentable categories is well $\lambda$-accessible.
This immediately follows from the fact that, for regular cardinals $\lambda\leq\mu$, an object of a locally 
$\lambda$-presentable category is $\mu$-presentable if and only if it is a $\mu$-small colimit of $\lambda$-presentable
objects. Recall that a category is $\mu$-\textit{small} if it has less than $\mu$ morphisms.

(3) Every finitely accessible functor is well finitely accessible. This follows from \cite{AR}, 2.13 (1) and the fact
that, for regular cardinals $\lambda\triangleleft\mu$, an object of a $\lambda$-presentable category is $\mu$-presentable 
if and only if it is a $\lambda$-directed $\mu$-small colimit of $\lambda$-presentable objects (see \cite{MP}, 2.3.11). 

(4) Let $I$ be an infinite set and consider the functor $\Set\rightarrow\Set$ that sends $X$ to $X^I$. This functor is accessible but not well 
accessible since there are arbitrarily large $|X|$ SUCH that $|X^I|>|X|$, hence the functor does not preserve $\mu$-presentable objects 
for all arbitrarily large enough $\mu$.  On the other hand, for $I$ finite, the functor is well accessible since it is finitely accessible. 
}
\end{rem}

A well accessible functor, by definition, will take $\mu$-presentable objects to $\mu$-presentable ones for all large enough $\mu$; however, it can 
lower presentability ranks.

\begin{exam}\label{ex3.3}
{
\em
Let $F:\Gr\to\Ab$ be the abelianization functor, i.e., the reflector from groups to abelian groups. Since $F$ is a left adjoint, it is well finitely 
accessible by \ref{re3.2} (2). There exist simple groups $G$ of arbitrarily large infinite cardinalities. Since $F(G)=0$ for such a $G$, $0$ is finitely 
presentable in $\Ab$ and $\rank(G)=|G|^+$, the functor $F$ can lower presentability ranks from $\kappa^+$ to $\omega$ for arbitrarily large $\kappa$.
}
\end{exam}

We say that a functor $F:\ck\to\cl$ \textit{reflects} $\lambda$-\textit{presentable objects} if $F(K)$ $\lambda$-presentable 
implies that $K$ is $\lambda$-presentable.

Recall that a morphism $g:B\to A$ is a split epimorphism if there exists $f:A\to B$ with $gf=\id_A$. Since, in this case, $g$
is a coequalizer of $fg$ and $\id_B$, $A$ is $\lambda$-presentable whenever $B$ is $\lambda$-presentable (see \cite{AR}, 1.16).

\begin{defi}\label{def3.4}
{
\em
We say that a functor $F:\ck\to\cl$ \textit{reflects split epimorphisms} if $f$ is a split epimorphism whenever $F(f)$ is a split
epimorphism.
 }
\end{defi} 

\begin{rem}\label{re3.5}
{
\em
Any functor $F:\ck\to\cl$ reflecting split epimorphisms is \textit{conservative}, i.e., it reflects isomorphisms. If all morphisms
of $\ck$ are monomorphisms then $F$ reflects split epimorphisms if and only if it is conservative.
}
\end{rem}

\begin{lemma}\label{le3.6}
Let $F:\ck\to\cl$ be a $\lambda$-accessible functor which reflects split epimorphisms. Then $F$ reflects $\lambda$-presentable 
objects. 
\end{lemma}
\begin{proof}
Let $F(K)$ be $\lambda$-presentable in $\cl$. Since $\ck$ is $\lambda$-accessible, $K$ is a $\lambda$-directed colimit
$(k_i:K_i\to K)_{i\in I}$ of $\lambda$-presentable objects. Since $F$ preserves $\lambda$-directed colimits and $F(K)$ 
is $\lambda$-presentable, there is $i\in I$ and $f:F(K)\to F(K_i)$ such that $F(k_i)f=\id_{F(K)}$. Since $F$ reflects split 
epimorphisms, this $k_i:K_i\to K$ is a split epimorphism. Thus $K$ is $\lambda$-presentable.
\end{proof}

\begin{propo}\label{prop3.7}
Let $F:\ck\to\cl$ be a well $\lambda$-accessible functor which reflects split epimorphisms. Then $F$ preserves presentability 
ranks $\mu$ for $\lambda<\mu$.
\end{propo}
\begin{proof}
Let $\rank(K)=\mu$, $\lambda<\mu$. Then $F(K)$ is $\mu$-presentable. Assume that $F(K)$ is $\nu$-presentable for some $\nu<\mu$. 
Without loss of generality, we can assume that $\lambda\leq\nu$. Following Lemma~\ref{le3.6}, $K$ is $\nu$-presentable, which 
is a contradiction. Thus $\rank(F(K))=\mu$.
\end{proof}

\section{Accessible categories with directed colimits}
An important class of accessible categories consists of accessible categories having directed colimits. It includes both $\infty,\omega$-elementary
categories and locally presentable ones. On the other hand, the basic $L_{\omega_1\omega_1}$ theory $T$ of well-ordered sets 
has both $\Mod(T)$ and $\Emb(T)$ $\omega_1$-accessible without having directed colimits. 

\begin{propo}\label{prop4.1}
Any accessible category with directed colimits is well accessible. 
\end{propo}
\begin{proof}
Let $\ck$ be a $\lambda$-accessible category with directed colimits and consider a regular cardinal $\lambda<\mu$. Given 
an object $K$ of $\ck$, there is a $\lambda$-directed colimit $(a_i:A_i\to K)_{i\in I}$ of $\lambda$-presentable objects $A_i$. 
Let $\hat{I}$ be the poset of all directed subsets of $I$ of cardinalities less than $\mu$ (ordered by inclusion). Clearly, 
$\hat{I}$ is $\mu$-directed. For each $M\in\hat{I}$, let $B_M$ be a colimit of a subdiagram indexed by $M$. Then $B_M$ is 
$\mu$-presentable. Since every subset of $I$ having less than $\mu$ elements is contained in a directed subset of $I$ having
less than $\mu$ elements (cf. \cite{AR}, 2.11), $K$ is a $\mu$-directed colimit of $B_M$, $M\in\hat{I}$. Thus $\ck$ is
$\mu$-accessible.
\end{proof}

\begin{lemma}\label{le4.2}
Let $\ck$ be a $\lambda$-accessible category with directed colimits and $K$ an object of $\ck$ which is not $\lambda$-presentable. 
Then $\rank(K)$ is a successor cardinal.
\end{lemma}
\begin{proof}
Let $\lambda<\mu=\rank(K)$. Following the proof of \ref{prop4.1}, $K$ is a $\mu$-directed colimit of objects $K_i$ which are directed 
colimits of size $\nu_i<\mu$ of $\lambda$-presentable objects. Since $K$ is $\mu$-presentable, it is a retract of some $K_i$. If
$\nu_i<\lambda$ then $K_i$ is $\lambda$-presentable and thus $K$ is $\lambda$-presentable as well. Let $\lambda\leq\nu_i$. Then
$K_i$ is $\nu_i^+$-presentable and thus $K$ is $\nu_i^+$-presentable too. Hence $\mu=\nu_i^+$.
\end{proof}

\begin{propo}\label{prop4.3}
Let $\ck$ and $\cl$ be accessible categories with directed colimits and $F:\ck\to\cl$ a functor preserving directed colimits.
Then $F$ is well accessible. 
\end{propo}
\begin{proof}
There is a regular cardinal $\lambda$ such that both $\ck$ and $\cl$ are $\lambda$-accessible and $F$ preserves 
$\lambda$-presentable objects (see \cite{AR} 2.19). Consider a regular cardinal $\lambda<\mu$ and let $K$ be a $\mu$-presentable 
object of $\ck$. Following the proof of \ref{prop4.1}, $K$ is a $\mu$-directed colimit of objects $B_M$ where each $B_M$ is 
a directed colimit of less than $\mu$ $\lambda$-presentable objects. Since $K$ is $\mu$-presentable, it is a retract of some 
$B_M$. Since $F(B_M)$ is a directed colimit of less than $\mu$ $\lambda$-presentable objects, $F(B_M)$ is $\mu$-presentable 
in $\cl$. Thus $F(K)$ is $\mu$-presentable in $\cl$ as a retract of $F(B_M)$. We have proved that $F$ is well $\lambda$-accessible.
\end{proof}

\begin{rem}\label{re4.4}
{
\em
Let $\Lin$ be the category of linearly ordered sets and order preserving injective mappings (they coincide with substructure embeddings). 
It is a finitely accessible category which is a ``minimal" $\infty,\omega$-elementary category in the sense that for every large 
$\infty,\omega$-elementary category $\ck$ there is a faithful functor $E:\Lin\to\ck$ preserving directed colimits; see \cite{MP}, 3.4.1. 
Its construction is based on Ehrenfeucht-Mostowski models. Following \ref{prop4.1} and \ref{prop4.3}, $E$ is well accessible. Since it is faithful, 
it reflects epimorphisms. Epimorphisms in $\Lin$ are isomorphisms and thus $E$ reflects split epimorphisms. Following \ref{prop3.7}, $E$ preserves 
presentability ranks $\mu$ starting from some cardinal $\lambda$. 
}
\end{rem}

Recall that a full subcategory $\ck$ of a category $\cl$ is called \textit{accessibly embedded} if there is a regular cardinal $\lambda$ 
such that $\ck$ is closed under $\lambda$-directed colimits in $\cl$. 

\begin{theo}\label{th4.5}
Accessible categories with directed colimits are precisely reflective and accessibly embedded subcategories of finitely accessible
categories.
\end{theo}
\begin{proof}
Every reflective and accessibly embedded subcategory of a fi\-ni\-te\-ly accessible category $\cl$ is accessible (see \cite{AR}, 
2.53) and has directed colimits calculated as reflections of directed colimits in $\cl$. Conversely, let $\ck$ be 
a $\lambda$-accessible category with directed colimits. Then $\ck$ is equivalent to a full subcategory $\Ind_{\lambda}(\cc)$ 
of $\Set^{\cc^{\op}}$ consisting of all $\lambda$-directed colimits of hom-functors (see \cite{AR}, 2.26); $\cc$ is a small 
category. Let $\Ind(\cc)$ be the full subcategory of $\Set^{\cc^{\op}}$ consisting of all directed colimits of hom-functors. 
Then $\Ind_{\lambda}(\cc)$ is closed in $\Ind(\cc)$ under $\lambda$-directed colimits. Given an object $X$ in $\Ind(\cc)$, 
we express it as a directed colimit of hom-functors and take their colimit $F(X)$ in $\Ind_{\lambda}(\cc)$ (recall that $\ck$ 
has directed colimits). Clearly, $F(X)$ is a reflection of $X$ in $\Ind_{\lambda}(\cc)$. Thus $\Ind_{\lambda}(\cc)$ is 
a reflective subcategory of $\Ind(\cc)$.
\end{proof}

\begin{theo}\label{th4.6}
Each large locally presentable category is LS-accessible. 
\end{theo}
\begin{proof}
Let $\ck$ be a locally presentable category. Then $\ck$ is locally $\lambda$-presentable for some regular cardinal $\lambda$
and, following \cite{AR} 5.30, $\ck$ is equivalent to $\Mod(T)$ for a limit theory of $L_{\lambda\lambda}(\Sigma)$ where $\Sigma$ 
is an $S$-sorted signature. Thus it suffices to prove that $\Mod(T)$ is LS-accessible.  

Let $\Set^S$ denote the category of $S$-sorted sets. The category $\Set^S$ is locally finitely presentable and, given an $S$-sorted
set $X=(X_s)_{s\in S}$, $\rank(X)=|X|^+$ where $|X|$ is defined as the cardinality of the disjoint union of $X_s$, $s\in S$. Let
$U:\ck\to\Set^S$ denote the forgetful functor. Since $\ck$ is large, there is $t\in S$ such that the sets $U(K)_t$ are arbitrarily large.
Let $V:\ck\to\Set$ be the composition of $U$ with the functor $\Set^S\to\Set$ sending $(X_s)$ to $X_t$. The functor $V$ preserves limits 
and $\lambda$-directed colimits (see the proof of \cite{AR} 5.9). By \cite{AR} 1.66, $V$ has a left adjoint $F$ which preserves 
$\mu$-presentable objects for each $\mu\geq\lambda$, cf.\ \cite{AR} Ex. 1.s(1). Thus $\rank(FX)\leq\rank(X)$ for each $X$ with $\lambda\leq\rank(X)$. 

Assume that $\Mod(T)$ is not LS-accessible. Then there is $X$ with $\lambda\leq\rank(X)$ and $\mu=\rank(FX)<\rank(X)$. Since $\Set^S$ 
is locally $\mu$-presentable, $X$ is a $\mu$-directed colimit of $\mu$-presentable objects $X_i$, $i\in I$. Let 
$$
(u_i:X_i\to X)_{i\in I}
$$ 
denote a colimit cocone. Since $F$ preserves colimits, 
$$
(Fu_i:FX_i\to FX)_{i\in I}
$$ 
is a colimit cocone. Since $FX$ is $\mu$-presentable, there is $j\in I$ and $r:FX\to FX_j$ such that $F(u_j)r=\id_{FX}$. Hence $Fu_j$ 
is a split epi\-mor\-phism. There is a mapping $f:X\to X$ such that $fu_j=u_j$ and $f\neq\id_X$. Since $F(f)UF(u_j)=F(u_j)$,
we have $F(f)=\id_{FX}$.   

Let $\eta:\Id\to VF$ be the adjunction unit. Since $\eta_Xf=VF(f)\eta_X=\eta_X$, $\eta_X$ is not a monomorphism. 
Since the sets $VK$ are arbitrarily large, there is $K$ in $\ck$ such that $|X|<|VK|$. Thus there is a monomorphism
$X\to VK$, which implies that $\eta_X$ is a monomorphism. We get a contradiction.
\end{proof}

\begin{rem}\label{re4.7}
{
\em
(1) We have proved the stronger result that the functor $F:\Set\to\ck$ preserves presentability ranks $\geq\lambda$.

(2) The argument even works for weakly locally presentable categories, i.e., for accessible categories with products. Such categories
are equivalent to $\Mod(T)$ where $T$ is a regular theory of $L_{\lambda\lambda}(\Sigma)$ (see \cite{AR} Ex. 5.e). The forgetful 
functor $U:\ck\to\Set^S$ has a weak left adjoint $F:Set^S\to\ck$ equipped with a natural transformation $\eta:\Id\to UF$. This follows
from the fact that $\ck$ is a small-injectivity class in a locally presentable category $\cl$ (see \cite{AR} 4.8) and thus it is
naturally weakly reflective in $\cl$ (cf. \cite{AHRT}). This suffices for the argument.

(3) Let $\ck$ be a large accessible category with coproducts. Assume that $\ck$ is $\lambda$-accessible and take the coproduct $K$
of a (representative) set of all $\lambda$-presentable objects. Then the functor 
$$U=\ck(K,-):\ck\to\Set
$$ 
has arbitrarily large values. $U$ has a left adjoint $F$ given by
$$
FX=\coprod\limits_XK
$$
By applying the proof of \ref{th4.6} to $F$, we get that $\ck$ is LS-accessible.  
}
\end{rem}

\begin{exam}\label{ex4.8 }
{
\em
Consider a one-sorted signature $\Sigma$ given by a sort $s$ and an $\omega$-ary function symbol $f:s^\omega\to s$. Let $T$ be 
the $L_{\omega_1\omega_1}(\Sigma)$ theory saying that $f$ is a bijection.  Thus $T$ consists of the formula
$$
(\forall y)(\exists !x)(f(x)=y)
$$
where $x$ is an $\omega$-string of variables, and $\exists !$ denotes unique existence. By general facts 
about categories of models of $\exists !$-sentences (called \textit{li\-mit theories}), the category $\Mod(T)$ is locally 
$\omega_1$-presentable; see \cite{AR}, 5.30. Hence, in particular, it is accessible and has directed colimits.
 
If $\mu$ has cofinality $\omega$ then 
$$
\mu^\omega > \mu
$$
cf.\ \cite{J}, Corollary 4 of Theorem 17.  Thus $T$ does not have models of cardinalities of cofinality $\omega$. Since
$\mu^\omega=\mu$ whenever the cardinal $\mu$ is of the form $\nu^\omega$, $\Mod(T)$ is a large category. Thus there is a proper
class of cardinalities in which $T$ has a model and a proper class of cardinalities in which $T$ does not have a model. But,
following \ref{th4.6}, $\Mod(T)$ is LS-accessible. The point is that presentation ranks and cardinalities will never start
to coincide. 

One can change the theory $T$ to an equational $L_{\omega_1\omega_1}$ theory $T'$ in the signature $\Sigma'$ consisting of $f$ and $\omega$ unary 
function symbols $g_i:s\to s$. Equations of $T'$ state that $f:X^\omega\to X$ and $\langle g_i\rangle_{i\in\omega}:X\to X^\omega$ are inverse maps 
(thus, bijections). The functor $V$ of the proof of \ref{th4.6} is the underlying set functor $\Mod(T')\to\Set$ and the functor $F$ is the free algebra 
functor $\Set\to\Mod(T')$. Following the proof of \ref{th4.6}, $\rank(FX)=|X|^+$ for any uncountable set $X$.  For free $T$-algebras $FX$, the difference 
between $\rank(FX)$ and the size of the underlying set $|VFX|$ can become arbitrarily large.  This clarifies the last sentence of the preceding
paragraph.
}
\end{exam}

\begin{defi}\label{def4.9}
{
\em
A functor $H:\ck\to\cl$ is \textit{full with respect to isomorphisms} if for any isomorphism $f:HA\to HB$ there is an isomorphism 
$\overline{f}:A\to B$ such that $H(\overline{f})=f$. 
}
\end{defi}

\begin{rem}\label{re4.10}
{
\em
Any functor $H$ full with respect to isomorphisms is \textit{essentially injective on objects} in the sense that $HA\cong HB$ implies that $A\cong B$.
Any faithful functor full with respect to isomorphisms is conservative, i.e.\ reflects isomorphisms.

A functor $H$ is called \textit{transportable} if for an isomorphism $f:HA\to B$ there is a unique isomorphism $\overline{f}:A\to\overline{B}$ 
such that $H(\overline{f})=f$ (this includes $H\overline{B}=B$).
}
\end{rem}

\begin{theo}\label{th4.11} 
Let $\ck$ be a large accessible category with directed colimits admitting a full with respect to isomorphisms and  faithful functor 
into a finitely accessible category preserving directed colimits. Then $\ck$ is LS-accessible.
\end{theo}
\begin{proof}
Following \ref{re2.3} (2), every finitely accessible category admits a finitely accessible full embedding into $\Str(\Sigma)$
for some finitary signature $\Sigma$. Thus we can assume that there is a full with respect to isomorphisms and faithful functor preserving 
directed colimits $H:\ck\to\Str(\Sigma)$. Consider the pullback
$$
\xymatrix@=4pc{
\ck \ar[r]^{H} & \Str(\Sigma) \\
\cl \ar [u]^{\overline{G}} \ar [r]_{\overline{H}} &
\Emb(\Sigma) \ar[u]_{G}
}
$$
where $G$ is the inclusion. Since $G$ is transportable, the pullback is equivalent to a pseudopullback (see \cite{MP} 5.1.1). Hence  
$\cl$ is accessible (see \cite{MP} 5.1.6) and clearly has directed colimits. The functors $\overline{H}$ and $\overline{G}$ preserve directed 
colimits, are faithful and full with respect to isomorphisms. In fact, for an isomorphism $f:\overline{G}L_1\to\overline{G}L_2$ in $\ck$, there 
is an isomorphism $\overline{f}:\overline{H}L_1\to \overline{H}L_2$ in $\Emb(\Sigma)$ with $G(\overline{f})=H(f)$. Thus there is an isomorphism
$g:L_1\to L_2$ with $\overline{G}(g)=f$. The argument for $\overline{H}$ is the same. Consequently, both $\overline{H}$ and $\overline{G}$ 
are essentially injective on objects (see Remark~\ref{re4.10}) and $\overline{G}$ is surjective on objects. Following \ref{prop4.1} and \ref{prop4.3}, 
there is a regular cardinal $\lambda$ such that both $\overline{G}$ and $\overline{H}$ are well $\lambda$-accessible. Since $\overline{G}$ is
faithful and full with respect to isomorphisms, it is conservative and thus reflects split epimorphisms (see \ref{re4.10} and \ref{re3.5}).
 
Now, following \ref{prop3.7}, $\overline{G}$ preserves presentability ranks $\nu$ for $\lambda<\nu$. Thus it suffices to prove
that $\cl$ is $\lambda$-LS-accessible. Consider a cardinal $\mu\geq\lambda$. Since $\ck$ is large and $\overline{G}$ is
essentially injective on objects, $\cl$ is large. Since $\overline{H}$ is essentially injective on objects, 
there is an object $L$ in $\cl$ such that $\mu\leq |\overline{H}L|$. Since $\cl$ is $\mu^+$-accessible, $L$ is a $\mu^+$-directed colimit 
of $\mu^+$-presentable objects $L_i$, $i\in I$. Let $(l_i:L_i\to L)_{i\in I}$ denote a colimit cocone. There is an embedding $f:A\to\overline{H}L$ 
with $|A|=\mu$. Since $A$ is $\mu^+$-presentable, $f=\overline{H}(l_j)g$ for some $g:A\to\overline{H}L_j$. Hence $\mu\leq|\overline{H}L_j|$.
Since $\overline{H}$ preserves presentability ranks $\nu$ for $\lambda<\nu$ (following \ref{prop3.7}), $\overline{H}L_j$ is
$\mu^+$-presentable, i.e., $|\overline{H}L_j|\leq\mu$. Hence $|\overline{H}L_j|=\mu$, i.e., $\rank(\overline{H}L_j)=\mu^+$. Therefore
$\rank(L_j)=\mu^+$. We have proved that $\cl$ is $\lambda$-LS-accessible.
\end{proof}

Using Ehrenfeucht-Mostowski models, we can prove a stronger result. (Note that any finitely accessible category has a faithful, directed colimit 
preserving functor into $\Set$ and thus any category $\ck$ from \ref{th4.11} has this property.)

\begin{theo}\label{th4.12}
Let $\ck$ be a large accessible category with directed colimits equipped with a faithful functor $H:\ck\to\Set$ preserving
directed colimits. Then $\ck$ is LS-accessible.
\end{theo}
\begin{proof}
Consider a pullback
$$
\xymatrix@=4pc{
\ck \ar[r]^{H} & \Set \\
\cl \ar [u]^{\overline{G}} \ar [r]_{\overline{H}} &
\Emb(\Set) \ar[u]_{G}
}
$$
analogous to that in the proof of \ref{th4.11}. Again, $\cl$ is a $\lambda$-accessible category with directed colimits. Following \ref{th4.5},
$\cl=\Ind_\lambda\cc$ is a full reflective subcategory of the finitely accessible category $\Ind\cc$. Let $F:\Ind\cc\to\cl$ be a left adjoint
to the inclusion $\cl\to\Ind\cc$. The composition $H\overline{G}F:\Ind\cc\to\Set$ preserves directed colimits and thus it is uniquely
determined by its domain restriction on $\cc$. Since $F$ is the identity functor on its domain restriction on $\cc$, the domain restriction 
of $H\overline{G}F$ on $\cc$ is faithful. We will prove that $H\overline{G}F$ is faithful. Consider two distinct morphisms $f,g:X\to Y$
in $\Ind\cc$. There are objects $C,D\in\cc$ and morphisms $u:C\to X$, $v:D\to Y$ and $f',g':C\to D$ such that $fu=vf'$ and $gu=vg'$.
Hence $f'$ and $g'$ are distinct and thus $H\overline{G}F(f')$ and $H\overline{G}F(g')$ are distinct. Since 
$H\overline{G}F(v)=G\overline{H}F(v)$ is a monomorphism, $H\overline{G}F(fu)$ and $H\overline{G}F(gu)$ are distinct. Therefore 
$H\overline{G}F(f)$ and $H\overline{G}F(g)$ are distinct. We have proved that $H\overline{G}F$ is faithful and thus $\overline{G}F$ 
is faithful.

Since every finitely accessible category is $\infty,\omega$-elementary (see \cite{MP} 5.2.6), \ref{re4.4} provides a faithful functor
$E:\Lin\to\Ind\cc$ preserving directed colimits. Thus $\overline{G}FE$ is faithful and preserves directed colimits. Following \ref{re4.4},
this functor preserves presentability ranks starting from some cardinal. Thus $\cl$ is $LS$-accessible. In order to prove that $\ck$ is
LS-accessible, we have to show that $\overline{G}$ is conservative -- then one uses \ref{re3.5}, \ref{prop4.3} and \ref{prop3.7}. 
Assume that $\overline{G}(f)$ is an isomorphism. Then $H\overline{G}(f)$, and thus $\overline{H}(f)$ are isomorphisms. Hence $f$ is
an isomorphism. 
\end{proof}

\begin{rem}\label{re4.13}
{
\em
We have proved that each category $\ck$ from \ref{th4.12} admits a faithful functor $E:\cl\to\ck$ from a finitely accessible category $\cl$
which preserves directed colimits and is surjective on objects. Consequently, it admits a faithful functor $\Lin\to\ck$ preserving directed
colimits. Moreover, if all morphisms in $\ck$ are monomorphisms, we do not need a pullback from the proof of \ref{th4.12} and the functor $E$ is
then even surjective on morphisms.
}
\end{rem}

We will give an example of an accessible category $\cw$ with directed colimits with no faithful functor $\cw\to\Set$ preserving directed colimits.  
This shows that the two parts of the assumption on the category $\ck$ in \ref{th4.12} are independent.  At the same time, this category $\cw$ is 
an example of an accessible category with directed colimits that is not $\infty,\omega$ elementary. On the other hand, we do not know whether 
the assumptions of \ref{th4.12} are necessary for the conclusion. Possibly every accessible category with directed colimits is LS-accessible; 
possibly every accessible category is LS-accessible.

\begin{exam}\label{ex4.14}
{
\em
Let $\cw$ be the category of well-ordered sets where morphisms are either order preserving injective mappings or constant mappings. Since the category 
$\cw_0$ of well-ordered sets and substructure embeddings is $\omega_1$-accessible, $\cw$ is $\omega_1$-accessible as well. The category $\cw_0$ 
does not have all directed colimits.  For example, a countable chain 
$$
A_0\to A_1\to\dots A_n\to\dots
$$
such that each $A_{n+1}$, $0\leq n$ contains an element $a$ smaller than any $x\in A_n$, does not have a colimit in $\cw_0$. In fact, it
does not have any compatible cocone $A_n\to A$, $0\leq n$. In $\cw$, it does have a colimit: the one-element chain $1$.  One can show that $\cw$ 
has all directed colimits.

Assume that there exists a faithful functor $H:\cw\to\Set$ preserving directed colimits. Following \ref{re4.13}, there is a faithful functor
$\Lin\to\cw$. This is impossible because $\cw$ is automorphism rigid -- the only isomorphisms $X\to X$ in $\cw$ are identities.
}
\end{exam}

\section {Abstract elementary classes}
 
Consider the following hierarchy of accessible categories with directed colimits:
\begin{enumerate}
\item[(1)] finitely accessible categories,
\item[(2)] $\infty,\omega$-elementary categories,
\item[(3)] accessible categories with directed colimits admitting a faithful functor preserving directed colimits
           into a finitely accessible category,
\item[(4)] accessible categories with directed colimits.  
\end{enumerate}
Each class is contained in the next.
The original definition of class (2), $\infty,\omega$-elementary categories, is model-theoretic, as categories equivalent to the category of models 
and homomorphisms of theories in $L^*_{\infty,\omega}$.  There exist categorical descriptions too: these are the categories equivalent to the category 
of points of some Grothendieck topos; and, by \cite{MP} 5.2.6, exactly the 2-categorical limits of finitely accessible categories and finitely accessible
functors. Since any finitely accessible category is equipped with a faithful functor into $\Set$ preserving directed colimits, class (3) coincides 
with the categories from Theorem \ref{th4.12}. Following \ref{ex4.14}, the inclusion of (3) in (4) is proper. The inclusion of (1) in (2) is also proper. 
We are not aware of any example of a category belonging to (3) but not to (2). We do know that any category from (2) admits a full embedding preserving 
directed colimits into a finitely accessible category.

We will show that the categories of models-and-strong-embeddings coming from Shelah's abstract elementary classes are between (2) and (4).
We will also introduce a class of abstract elementary categories as ``abstract elementary classes'' without any assumption about morphisms 
being embeddings.

\begin{defi}\label{def5.1}
{
\em
A functor $H:\ck\to\cl$ will be called \textit{nearly full} if for each commutative triangle 
$$
\xymatrix@=3pc{
HA \ar[rr]^{H(h)}
\ar[dr]_{f} && HC\\
& HB \ar[ur]_{H(g)}
}
$$
there is $\overline{f}:A\to B$ in $\ck$ such that $H(\overline{f})=f$.

A subcategory $\ck$ of a category $\cl$ will be called \textit{nearly full} if the embedding $\ck\to\cl$ is nearly full. 
}
\end{defi} 

\begin{rem}\label{re5.2}
{
\em
(1) Every full functor is nearly full. If $H$ is faithful and nearly full, we also have $g\overline{f}=h$ in the definition above.
Thus every faithful and nearly full functor reflects split epimorphisms.

(2) \cite{Ki} calls nearly full subcategories coherent.  

(3) Any subcategory $\ck$ of a finitely accessible category $\cl$ closed under directed colimits is \textit{replete} and full with respect 
to isomorphisms. In fact, being closed under directed colimits means that each directed colimit in $\cl$ of objects from $\ck$ belongs 
to $\ck$. And, directed colimits are determined up to an isomorphism. Recall that replete means to be closed under isomorphic objects.

If $\ck$ is a subcategory of an accessible category $\cl$ and the embedding of $\ck$ to $\cl$ preserves directed colimits then 
the \textit{replete closure} $\overline{\ck}$ of $\ck$ is closed under directed colimits in $\cl$.  Moreover the categories $\ck$
and $\overline{\ck}$ are equivalent.
}
\end{rem}

\begin{defi}\label{def5.3}
{
\em
An accessible categories with directed colimits will be called an \textit{abstract elementary category} if it admits a full with respect 
to isomorphisms and nearly full embedding preserving directed colimits into a finitely accessible category.
}
\end{defi}
 
Abstract elementary categories are closely related to Shelah's Abstract Elementary Classes; they differ from those introduced in \cite{Ki}.
We recall the definition of an abstract elementary class using the language of category theory. 

\begin{defi}\label{def5.4}
{
\em
Let $\Sigma$ be a finitary signature. A nearly full subcategory $\ck$ of $\Emb(\Sigma)$ is called an \textit{abstract elementary class} 
if it is closed in $\Emb(\Sigma)$ under directed colimits and there is a cardinal $\lambda$ such that if $f:A\to B$ is a substructure 
with $B\in\ck$ then there is $h:A'\to B$ in $\ck$ such that $f$ factorizes through $h$ and $|A'|\leq |A|+\lambda$.
}
\end{defi}

The standard formulation (see \cite{B}) uses colimits of continuous chains instead of directed colimits. But it is well known that 
this does not change anything (see \cite{B}, or \cite{AR} 1.7; continuous chains are called smooth there). The standard formulation 
also includes that $\ck$ is replete and full with respect to isomorphisms but \ref{re5.2} (3) shows that our formulation is the same. 
We allow many-sorted signatures in \ref{def5.4}. Since each many-sorted signature can be made single-sorted, this does not change the concept 
of an abstract elementary class. 

\begin{theo}\label{th5.5}
Let $\Sigma$ be a finitary signature. A nearly full subcategory $\ck$ of $\Emb(\Sigma)$ is an abstract elementary class 
if and only if it is an accessible category closed under directed colimits in $\Emb(\Sigma)$. 
\end{theo}
\begin{proof}
Let $\ck$ be an abstract elementary class in $\Emb(\Sigma)$. We know that $\Emb(\Sigma)$ is finitely accessible and there is a regular 
cardinal $\lambda<\kappa$ such that, for each regular cardinal $\kappa\leq\mu$, $\mu$-presentable objects in $\Emb(\Sigma)$ are precisely 
$\Sigma$-structures $A$ such that $|A|<\mu$. Since $\Emb(\Sigma)$ is $\kappa$-accessible, Definition \ref{def5.4} yields that each object 
of $\ck$ is a $\kappa$-directed colimit of $\kappa$-presentable objects in $\ck$. Thus $\ck$ is $\kappa$-accessible.

Conversely, let $\ck$ be an accessible nearly full subcategory of $\Emb(\Sigma)$ closed under directed colimits. Then the embedding
$\ck\to\Emb(\Sigma)$ is well $\lambda$-accessible for some regular cardinal $\lambda$ greater than $\kappa$ above. Let $f:A\to B$ be 
a substructure embedding with $B\in\ck$ and put $\mu=|A|+\lambda$. Since $\ck$ is $\mu^+$-accessible, $B$ is a $\mu^+$-directed colimit of 
$\mu^+$-presentable objects $B_i$, $i\in I$ in $\ck$. Since the embedding $\ck\to\Emb(\Sigma)$ preserves $\mu^+$-directed colimits
and $\mu^+$-presentable objects, $B$ is a $\mu^+$-directed colimit of objects $B_i$ which are $\mu^+$-presentable in $\Emb(\Sigma)$.  
Since $A$ is $\mu^+$-presentable (because $\kappa<\mu^+$), the substructure embedding $f:A\to B$ factorizes through some $B_j$, $j\in I$. It
suffices to put $A'=B_j$.
\end{proof}

This result is an improvement of \cite{L} 5.9 and 5.10 (see \cite{L1} 4.1 and 4.9, as well; there is a related work \cite{Ki}). 
Lieberman assumes that $\ck$ is well accessible while we use \ref{prop4.1} to prove this from the existence of directed colimits.
Otherwise, our proof is the same as that of Lieberman. 
 
\begin{rem}\label{re5.6}
{
\em
(1) Categories $\Elem(T)$ where $T$ is an $L_{\kappa\omega}$ theory and $\Emb(T)$ where $T$ is a basic $L_{\kappa\omega}$ theory
are abstract elementary classes.

(2) Each large abstract elementary class is an LS-accessible category. This follows from \ref{th4.11} but also directly from 
the proof of \ref{th5.5}. Following this proof, we have $\rank(A)=|A|^+$ for each $A$ in $\ck$ with $|A|\geq\lambda$. 
}
\end{rem}

We can write our characterization of abstract elementary classes in the language of category theory, i.e., without using 
$\Sigma$-structures.

\begin{coro}\label{cor5.7} 
A category is equivalent to an abstract elementary class if and only if it is an accessible category with directed colimits whose morphisms 
are monomorphisms and which admits a full with respect to isomorphisms and nearly full embedding into a finitely accessible category preserving 
directed colimits and monomorphisms.
\end{coro}
\begin{proof}
Necessity is evident because $\Emb(\Sigma)$ is finitely accessible. Consider a category $\ck$ satisfying the conditions above. Let
$H:\ck\to\cl$ be the corresponding functor into a finitely accessible category. Consider the canonical functor
$$
E:\cl\to\Set^{\ca^{\op}}.
$$
(see \cite{AR}, 1.25). $E$ preserves directed colimits (see \cite{AR}, 1.26). Since objects of $\Set^{\ca^{\op}}$ can be viewed as many-sorted 
unary algebras, $\Set^{\ca^{\op}}$ is a full subcategory of $\Str(\Sigma)$ for a finitary signature $\Sigma$ containing operation symbols only. 
Thus embeddings in $\Str(\Sigma)$ coincide with monomorphisms. Hence the codomain restriction of the composition $EH$ is a full with respect 
to isomorphisms and nearly full embedding  
$$
\ck\to\Emb(\Sigma)
$$
preserving directed colimits. Following \ref{re5.2} (3) and \ref{th5.5}, $\ck$ is equivalent to an abstract elementary class.  
\end{proof}

\begin{rem}\label{re5.8}
{
\em
(1) Our definition of abstract elementary category is motivated by Corollary 5.7; it results from dropping the hypotheses on monomorphisms.
Any finitely accessible category whose morphisms are not monomorphisms, like posets and isotone mappings, is example of an abstract elementary category 
that is not an abstract elementary class. 

(2) Any abstract elementary class is an abstract elementary category. On the other hand, let $\ck$ be an abstract elementary category with
a functor $H$ into a finitely accessible category from \ref{def5.3}. Without any loss of generality, we can assume that $H:\ck\to\Str(\Sigma)$
where $\Sigma$ is a finitary signature. Consider the pullback from the proof of \ref{th4.11}. Then $\cl$ is an accessible category with directed 
colimits whose morphisms are monomorphisms and $\overline{H}$ is a full with respect to isomorphisms and nearly full embedding preserving directed 
colimits and monomorphisms. Thus $\cl$ is equivalent to an abstract elementary class. Moreover, the functor $\overline{G}$ is onto on objects 
and preserves presentability ranks starting from some regular cardinal.

(3) We are not aware of any abstract elementary category which is not $\infty,\omega$-elementary. A simple example of an abstract 
elementary class $\ck$ in $\Sigma$ which is not closed under $L_{\infty,\omega}$-elementary equivalence (thus not axiomatizable
in any $L_{\kappa\omega}(\Sigma)$) is given in \cite{K}, 2.10. There, $\Sigma$ is a single-sorted signature containing just a unary 
relation symbol $P$. Objects of $\ck$ are $\Sigma$-structures $K$ such that $P_K$ is countable and the complement of $P_K$ is infinite.
Morphisms of $\ck$ are substructure embeddings which are identities on $P$. But, $\ck$ is isomorphic to the category of infinite sets and monomorphisms. 
The latter category is axiomatizable by a basic $L_{\omega\omega}$ theory in the empty signature. Thus it is $\infty,\omega$-elementary.
}
\end{rem}

\begin{propo}\label{prop5.9}
Let $\ck$ be an abstract elementary category. Then there is an $\infty,\omega$-elementary category $\cl$ and a faithful functor $E:\cl\to\ck$ 
which preserves directed colimits and is surjective on objects.
\end{propo}
\begin{proof}
Let $\overline{G}:\cl\to\ck$ be the embedding of an abstract elementary class $\cl$ from \ref{re5.8} (2). Shelah's Presentation Theorem
(see \cite{B} 4.15) yields a finitary signature $\Sigma\subseteq\Sigma'$, an $L_{\omega\omega}$-theory $T'$  and a set $\Gamma$ 
of quantifier free types in $\Sigma'$ such that $\cl$ consists of $\Sigma$-reducts of $T'$-models omitting types from $\Gamma$. Omitting 
a type can be expressed as an $L_{\infty\omega}$-sentence. By adding these sentences for types from $\Gamma$ to $T'$, we get 
an $L_{\infty\omega}$-theory $T''$. Since $\Elem(T'')$ is $\infty,\omega$-elementary (see \ref{re2.3} (5), the reduct functor $R:\Elem(T'')\to\cl$ 
has the desired properties for $\cl$. Thus $\overline{G}R$ has these properties for $\ck$.
\end{proof}

\section{Categoricity}
\begin{defi}\label{def6.1}
{
\em
Let $\lambda$ be an infinite cardinal. A category $\ck$ is called $\lambda$-\textit{categorical} it it has,
up to isomorphism, precisely one object of the presentability rank $\lambda^+$.
}
\end{defi} 

\begin{rem}\label{re6.2}
{
\em
Following \ref{re5.6} (2), the definition \ref{def6.1} (suggested in \cite{R2}) is in accordance with its model theoretic meaning 
in abstract elementary classes for sufficiently large cardinals $\lambda$.
}
\end{rem}

Shelah's Categoricity Conjecture claims that for every abstract elementary class $\ck$ there is a cardinal $\kappa$ such that $\ck$ 
is either $\lambda$-categorical for all cardinals $\kappa\leq\lambda$ or $\ck$ is not $\lambda$-categorical for any cardinal $\kappa\leq\lambda$.
Following \ref{re6.2}, this is a statement about $\ck$ as a category, i.e., it does not depend on the signature in which the abstract elementary class 
is presented. 

By Remark \ref{re5.8} (2), Shelah's Categoricity Conjecture is equivalent to the Categoricity Conjecture for abstract elementary categories.  
It is natural to ask about the status of the Categoricity Conjecture for accessible categories belonging to other levels of our hierarchy.  At present, 
we can offer little information other than this easy observation.
 
\begin{exam}\label{ex6.3}
{
\em
Suppose $\ck$ is an accessible category which is not LS-acce\-ssib\-le.  Let $\cl=\ck\sqcup\Set$ be the disjoint union of $\ck$ and the category 
of sets. Then $\cl$ is an accessible category and there is a proper class of cardinals $\lambda$ such that $\cl$ is $\lambda$-categorical and, 
at the same time, a proper class of cardinals $\lambda$ such that $\cl$ is not $\lambda$-categorical.

If this $\ck$ has directed colimits then the categoricity conjecture fails for the class (4).  By \ref{th4.12}, however, this simple trick does not help for categories in class (3) of our hierarchy.
}
\end{exam}

\begin{rem}\label{re6.4}
{
\em
Shelah's conjecture seems to be unknown even for finitely accessible categories, lying at level (1) of the hierarchy. It would be interesting 
to understand whether the exquisite Galois-theoretic machinery of \cite{B} can be brought to bear implications in this setting. 
}
\end{rem}

\end{document}